\documentclass[a4paper,12pt]{amsart}

\usepackage{fouriernc}
\usepackage[T1]{fontenc}

\usepackage{amssymb,latexsym,amsmath,epsfig,amsthm} 

\usepackage{bbm}
\usepackage{enumerate}
\usepackage{rotating}

\usepackage{geometry}
\newtheorem{proposition}{Proposition}

\theoremstyle{definition}

\def\CC{{\mathbb C}}

\def\MM{{\mathbb M}}
\def\NN{{\mathbb N}}
\def\QQ{{\mathbb Q}}

\def\RR{{\mathbb R}}

\def\ZZ{{\mathbb Z}}

\def\vecp{{\text{\boldmath$p$}}}

\def\vecx{{\text{\boldmath$x$}}}

\def\vecalf{{\text{\boldmath$\alpha$}}}

\def\vecomega{{\text{\boldmath$\omega$}}}

\def\vecnull{{\text{\boldmath$0$}}}

\def\scrB{{\mathcal B}}

\def\scrD{{\mathcal D}}

\def\fB{{\mathfrak B}}

\def\Re{\operatorname{Re}}

\def\e{\mathrm{e}}

\def\SL{\operatorname{SL}}

\def\vol{\operatorname{vol}}

\def\GamG{\Gamma\backslash G}

\def\trans{\,^\mathrm{t}\!}

\numberwithin{equation}{section}

\begin{document}

\begin{center}
\uppercase{\bf \boldmath The log moments of smallest denominators}
\vskip 20pt
{\bf Jens Marklof\\
{\it School of Mathematics, University of Bristol, Bristol BS8 1UG, U.K.}\\
{\tt j.marklof@bristol.ac.uk}}
\end{center}
\vskip 20pt
\centerline{\it To appear in Integers} 
\vskip 30pt

\centerline{\sc Abstract}
\vskip 8pt
\noindent
This paper studies the logarithmic moments of the smallest denominator of all rationals in a shrinking interval with random center. Convergence follows from the more general results in [arXiv:2310.11251, Bull. Lond. Math. Soc., to appear], and the key point of this note is the derivation of explicit formulas for the moments of the limit distribution in dimension one. This answers questions raised by Meiss and Sander in their numerical study of minimal resonance orders for torus maps with random rotation vectors [arXiv:2310.11600].

\thispagestyle{empty}
\baselineskip=12.875pt
\vskip 30pt

\section{Introduction}\label{secIntro}

Consider the smallest denominator of all fractions in an interval of length $\delta$ centered at $x$,  
\begin{equation}\label{oneptone}
q_{\min}(x,\delta) = \min \left\{ q\in\NN : \exists \tfrac{p}{q} \in\QQ \cap (x-\tfrac{\delta}{2},x+\tfrac{\delta}{2}) \right\}.
\end{equation}
In their investigation of the breakdown of invariant tori in integrable systems, Meiss and Sander \cite{Meiss21,Meiss23} carried out a numerical study of the distribution of $\log q_{\min}(x,\delta)$ for random $x$ (as well as higher dimensional variants, which we return to in Section \ref{secMS}), and asked for a proof of a limit law and specifically the convergence of its expectation value as $\delta\to 0$. The asymptotics of the expectation value of $q_{\min}(x,\delta)$ (without taking the log) was already known due to work of Chen and Haynes \cite{Chen23}.  

We start with the following limit law, which is equivalent to \cite[Proposition 1]{Smallest} (replace $L$ by $\e^L$ and take $\eta_{\log}(s)=\e^s\eta(\e^s)=\tfrac{6}{\pi^2}\, \e^{2s}\, H(\tfrac{3}{\pi^2} \e^{2s})$, where $\eta(s)$ is the limiting density of smallest denominators and $H(s)$ is the Hall distribution as defined in \cite{Smallest}):

\begin{proposition}
For any interval $\scrD\subset[0,1]$ and $L\in\RR$, we have
\begin{equation}\label{LD1}
\lim_{\delta\to 0} \vol\left\{ x\in \scrD :   \log q_{\min}(x,\delta) > L - \tfrac12\log \delta  \right\} = \vol\scrD \; \int_L^\infty \eta_{\log}(s)\, ds
\end{equation}
with the probability density
\begin{equation}\label{formula2}
\eta_{\log}(s) = \tfrac{6}{\pi^2}
\begin{cases}
\e^{2s}  & \text{if $s\leq 0$}\\
-\e^{2s}+2+4 s & \text{if $0\leq s \leq \log 2$}\\
-\e^{2s}+ 2+2 \e^{2s} \sqrt{\frac14 - \e^{-2s}}-4\log\bigg(\frac12+\sqrt{\frac14-\e^{-2s}}\bigg) & \text{if $s\geq \log 2$.}
\end{cases}
\end{equation}
\end{proposition}

\begin{figure}
\begin{center}
\includegraphics[width=0.7\textwidth]{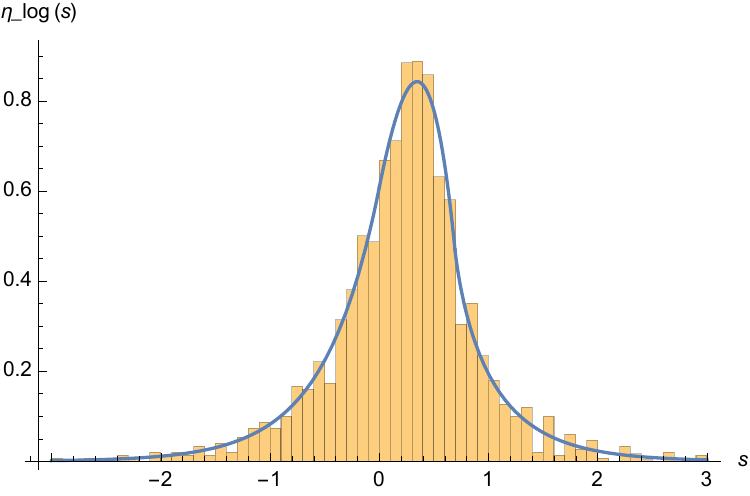}
\end{center}
\caption{The limit density $\eta_{\log}(s)$ compared to the distribution of the logarithm of the smallest denominator of rationals in each interval $[\frac{j}{3000},\frac{j+1}{3000})$, $j=0,\ldots,2999$.} \label{fig1}
\end{figure}
Note that the mode of $\eta_{\log}(s)$ is $s=\frac12\log 2$ (cf.~Figure~\ref{fig1}), and its tails are
\begin{equation}
\eta_{\log}(s) 
\begin{cases}
=\frac{6}{\pi^2} \e^{-2|s|} & (s\leq 0)\\
\sim\frac{12}{\pi^2} \e^{-2s} & (s\to\infty) .
\end{cases}
\end{equation}
This follows from the tail estimates of the Hall distribution; see \cite{Smallest} for details.

We now turn to the convergence of logarithmic moments, which follows directly from the convergence of moments proved in \cite{Smallest}.

\begin{proposition}
For any interval $\scrD\subset[0,1]$ and $n\in\ZZ_{\geq 0}$, we have
\begin{equation}\label{LD1Mom}
\lim_{\delta\to 0}   \int_{\scrD}   \left(\log q_{\min}(x,\delta)+\tfrac12\log\delta\right)^n dx=
\mu_n \vol\scrD , 
\end{equation}
with
\begin{equation}\label{LD1Mom222}
\mu_n = \int_0^\infty s^n \eta_{\log}(s)\, ds = \int_0^\infty (\log s)^n \eta(s)\, ds .
\end{equation}
\end{proposition}

\begin{proof}
This follows by dominated convergence from \cite[Proposition 2]{Smallest}, as the moments of the logarithm are bounded above by the first positive plus first negative moment of smallest denominators, both of which have a finite limit. 
\end{proof}

The above limit theorems for the distribution and logarithmic moments also hold (with identical limits) when $x$ is sampled over the discrete set $x=x_0+\frac{j}{N}$ ($j=1,\ldots,N$) with $x_0$ fixed, provided $N^{-1} = O(\delta)$ as $\delta\to 0$ and $N\to\infty$; see Figure~\ref{fig1}. This follows by the same argument as in the continuous sampling case, using now \cite[Propositions 5 and 6]{Smallest}. We refer the reader to \cite{Balazard23,Smallest,Shparlinski23} for more background and results in this setting.

In higher dimensions, the proof of convergence of logarithmic moments of smallest denominators of rational vectors follows similarly from the results in \cite[Section 2]{Smallest}. 
Section \ref{secMS} of this note provides the asymptotics for Meiss and Sander's minimal resonance orders \cite{Meiss21,Meiss23}, a different higher dimensional variant of smallest denominators.

The main point of the present paper is the calculation of explicit formulas for the moments $\mu_n$. We first note that the moment generating function of the limit distribution $\eta_{\log}(s)$ is, for $|\Re\alpha|<2$, $\alpha\neq 0$,
\begin{equation}
M(\alpha)  = \frac{24}{\pi^2 \alpha (\alpha+2)} \left( \frac{2}{\alpha} + 2^{\alpha} \mathrm{B}\left(-\frac{\alpha}{2},\frac12\right)\right)  
\end{equation}
where $\mathrm{B}(x,y)$ is the beta function (Euler's integral of the first kind). Recall from \cite{Smallest} that $M(\alpha)$ represents the  (complex) $\alpha$-moment of the density of the limit distribution $\eta(s)$ of the small denominators, and are closely related to the moments of the distance function for the Farey sequence determined by Kargaev and Zhigljavsky \cite{Kargaev97}.
Indeed, $M(\alpha)$ has an analytic continuation to $\alpha=0$ and we have, for $n=1,2,3,\ldots$,
\begin{equation}
\mu_n  = \left. \frac{d^n}{d\alpha^n} \int_0^\infty s^\alpha \eta(s)\, ds  \right|_{\alpha=0}
=  \left. \frac{d^n M(\alpha)}{d\alpha^n} \right|_{\alpha=0}.
\end{equation}
Furthermore, \cite[Proposition 2]{Smallest} is synonymous with the convergence of the moment generating function for the logs. In Section \ref{secExplicit} we prove the following explicit formula for $\mu_n$,
\begin{equation}\label{inv0009}
\mu_n = \frac{n!}{2}  \sum_{k=0}^n \sum_{j=0}^k   \frac{(\log 2)^{n-k}}{(n-k)!}\;  \frac{(-1)^{j} }{2^j(k-j)!} \; \rho_{k-j}  ,
\end{equation}
where
\begin{equation}\label{rhon}
\rho_n = \frac{24}{\pi^2} \left( \frac{n!}{2^{n+1} } \sum_{j=1}^\infty \begin{pmatrix} 2j-1 \\ j\end{pmatrix} \frac{1}{2^{2j-1} j^{n+2}}  + \frac{2(-\log 2)^{n+2}}{(n+2)(n+1)}  \right) .
\end{equation}
For $n=1,2,3$ this specializes to the surprisingly simple 
\begin{equation}\label{explicit}
\begin{split}
\mu_1 & =\frac{6}{\pi^2} \zeta (3) -\frac12=0.230763\ldots \\
\mu_2 & = \frac{3\pi^2}{40}+\frac12-\frac{6}{\pi^2} \zeta (3) =0.509457\ldots \\
\mu_3 & = \left( \frac{9}{\pi^2} -\frac32\right) \zeta (3)+\frac{27}{\pi^2} \zeta (5)-\frac34 -\frac{9 \pi ^2}{80} = 0.269423\ldots
\end{split}
\end{equation}
where $\zeta(s)$ is the Riemann zeta function. The proof of \eqref{explicit} uses a recursion formula (cf.~\eqref{bbcc}, Section \ref{secExplicit}) rather than a direct evaluation of \eqref{rhon}. 
We thus have for the standard deviation
\begin{equation}\label{explicit2}
\sigma = \left(\mu_2 - \mu_1^2\right)^{1/2} = \left(\frac{3 \pi ^2}{40}+\frac{1}{4}  -\frac{36}{\pi ^4} \zeta (3)^2 \right)^{1/2} = 0.67543\ldots
\end{equation}
and skewness
\begin{equation}\label{explicit3}
\begin{split}
\gamma & = \frac{\mu_3 - 3 \mu_1 \mu_2 + 2\mu_1^3}{\sigma^3} \\ 
& = 
\sigma^{-3} \left(\frac{432}{\pi ^6} \zeta (3)^3-\frac{57}{20}  \zeta (3)+\frac{27}{\pi ^2} \zeta (5)-\frac{1}{4}\right) \\
& = -0.190475\ldots .
\end{split}
\end{equation}

Meiss and Sander \cite{Meiss23} used $\log_{10}$ rather than $\log$, and also intervals of length $2\delta$ instead of $\delta$. Hence $\delta^{1/2} q_{\min}^\text{MS}(x,\delta)  = 2^{-1/2} (2\delta)^{1/2} q_{\min}(x,2\delta)$, and the necessary adjustments lead to 
\begin{equation}
\mu_1^\text{MS}=\frac{\mu_1-\frac12 \log 2}{\log 10} = -0.0502959\ldots, 
\end{equation}
\begin{equation}
\sigma^\text{MS}=\frac{(\mu_2-\mu_1^2)^{1/2}}{\log 10} = 0.293336\ldots,
\end{equation}
which is compatible with the numerical results of \cite{Meiss23}, $\mu_1^\text{MS}=-0.05 \pm 0.001$ and $\sigma^\text{MS}=0.2935\pm0.0006$.

As noted in \cite{Smallest}, the distribution of smallest denominators is closely related to the void distribution for the Farey sequence \cite{Boca05,Kargaev97}, as well as to the directional statistics of Euclidean lattice points \cite{Boca06,partI}. All three have the same limiting statistics, so our formulas for the moments apply in these settings. Similar limit distributions also arise in the study of the free path length in the periodic Lorentz gas. The connection between the two is explained in \cite{partI}, and indeed the entropy formulas established by Boca and Zaharescu \cite{Boca07} are similar to $\mu_1$, including the appearance of $\zeta(3)$.

\section{Explicit Formulas for Logarithmic Moments}\label{secExplicit}

To calculate the limiting logarithmic moments $\mu_n$ in dimension one, it will be convenient to use a slightly different normalization, shifting the log of the smallest denominator by $\log 2$, so that
\begin{equation}\label{LD1Mom828}
\lim_{\delta\to 0}   \int_{\scrD}   \left(\log q_{\min}(x,\delta)+\tfrac12\log\delta-\log 2\right)^n dx=
\tilde\mu_n \vol\scrD .
\end{equation}
The above convergence is implied by \eqref{LD1Mom} with 
\begin{equation}
\tilde\mu_n = \int_0^\infty (\log s-\log 2)^n \eta(s)\, ds =\sum_{k=0}^n \begin{pmatrix} n \\ k \end{pmatrix} (-\log 2)^{n-k} \mu_k.
\end{equation}
We can invert this relation to get a formula for the original moments,
\begin{equation}\label{inv0}
\mu_n = \sum_{k=0}^n \begin{pmatrix} n \\ k \end{pmatrix} (\log 2)^{n-k} \tilde\mu_k.
\end{equation}
Now set
\begin{equation}\label{principe}
P_n(x) = 2 \left(x-\log 2\right)^n + n \left(x-\log 2\right)^{n-1} 
\end{equation}
and observe that
\begin{equation}
P_n(\log s) = \frac{d^n}{d\alpha^n} \left[ (\alpha+2)  \left(\frac{s}{2}\right)^\alpha \right]_{\alpha=0}
= 2 \left(\log\frac{s}{2}\right)^n + n \left(\log\frac{s}{2}\right)^{n-1} .
\end{equation}
By \eqref{LD1Mom} we therefore have
\begin{equation}\label{LD1Mom007}
\lim_{\delta\to 0}   \int_{\scrD}   P_n\left(\log q_{\min}(x,\delta)+\tfrac12\log\delta\right) dx=
\rho_n \vol\scrD , 
\end{equation}
where
\begin{equation}
\begin{split}
\rho_n & = \int_0^\infty P_n(\log s)\, \eta(s)\, ds \\
& = \left. \frac{d^n}{d\alpha^n} \int_0^\infty (\alpha+2)  \left(\frac{s}{2}\right)^\alpha \eta(s)\, ds  \right|_{\alpha=0} \\
& =  \left. \frac{d^n \MM(\alpha)}{d\alpha^n} \right|_{\alpha=0},
\end{split}
\end{equation}
and
\begin{equation}
\MM(\alpha) = (\alpha+2)  2^{-\alpha} M(\alpha) .
\end{equation}
Relation~\eqref{principe} implies 
\begin{equation}
\rho_n = 2 \tilde\mu_n + n \tilde\mu_{n-1}, \quad \rho_0=2 .
\end{equation}
This can be inverted,
\begin{equation}\label{inversion}
\tilde\mu_n = \frac{n!}{2} \sum_{j=0}^n (-2)^{-j}  \frac{\rho_{n-j}}{(n-j)!} .
\end{equation}
Relation~\eqref{inv0009} now follows from \eqref{inv0} and \eqref{inversion}.  

The next task is thus to work out $\rho_n$.
We organize the terms in $\MM(\alpha)$ as follows:
\begin{equation}\label{MMalpha}
\begin{split}
\MM(\alpha)  & = \frac{24}{\pi^2}\;\frac1\alpha \left( 2^{-\alpha} \frac{2}{\alpha} + \mathrm{B}\left(-\frac{\alpha}{2},\frac12\right)\right)  \\
& = \frac{24}{\pi^2}\left[ \frac1\alpha  \left(\frac{2}{\alpha} - 2\log 2 + \mathrm{B}\left(-\frac{\alpha}{2},\frac12\right)\right)  + 2\;\frac{2^{-\alpha} - 1 + \alpha \log 2}{\alpha^2} \right] .
\end{split}
\end{equation}
We will express $\rho_n$ in terms of the coefficients $b_k$ of the Laurent expansion of the beta function,
\begin{equation}\label{betaLaurent}
\mathrm{B}\left(z,\frac12\right)= \sum_{k=-1}^\infty b_k z^k ,
\end{equation}
that is, 
\begin{equation}\label{betaLaurent2}
b_{-1}=1,\qquad b_n =  \frac{1}{n!} \; \frac{d^n}{dz^n} \left[ \mathrm{B}\left(z,\frac12\right) - \frac{1}{z} \right]_{z=0} \quad (n\geq 0).
\end{equation}
Then
\begin{equation}\label{betaLaurent3}
\frac1\alpha  \left(\frac{2}{\alpha} - 2\log 2 + \mathrm{B}\left(-\frac{\alpha}{2},\frac12\right)\right)  
= \sum_{k=0}^\infty (-1)^{k+1} 2^{-k-1} b_{k+1} \alpha^k
\end{equation}
and hence
\begin{equation}\label{betaLaurent4}
 \frac{d^n}{d\alpha^n}  \left[\frac1\alpha  \left(\frac{2}{\alpha} - 2\log 2 + \mathrm{B}\left(-\frac{\alpha}{2},\frac12\right)\right)  \right]_{\alpha=0}
= n!  (-1)^{n+1} 2^{-n-1} b_{n+1} .
\end{equation}
We view \eqref{betaLaurent} and \eqref{betaLaurent3} here as a formal series (ignoring their convergence), and may take  \eqref{betaLaurent2} as the actual definition of the $b_n$ from which \eqref{betaLaurent4} follows.

In order to work out an explicit formula for $b_n$, we note that by \cite[8.382.3]{GR} we have
\begin{equation}\label{Bexp}
\mathrm{B}\left(z,\frac12\right)= \frac{1}{z}+\sum_{j=1}^\infty \frac{(2j-1)!!}{2^j j!} \; \frac{1}{z+j} ,
\end{equation}
where the series converges (in view of Sterling's formula) absolutely and uniformly on compacta excluding $z= -1,-2,-3,\ldots$.
Using the $n$th derivative of this expression as an input to \eqref{betaLaurent2} yields
\begin{equation}
b_n =  (-1)^n \sum_{j=1}^\infty \frac{(2j-1)!!}{2^j j! j^{n+1}} = (-1)^n \sum_{j=1}^\infty \begin{pmatrix} 2j-1 \\ j\end{pmatrix} \frac{1}{2^{2j-1} j^{n+1}} ,
\end{equation}
which converges for all $n\geq 0$, again by Sterling's formula.

We can alternatively calculate $b_n$ recursively with help from the formula
\begin{equation}\label{derry}
\begin{split}
\frac{d}{dz} \mathrm{B}\left(z,\frac12\right) 
&= \left(\psi(z)-\psi(z+\tfrac12)\right) \mathrm{B}\left(z,\frac12\right) \\
& = \left( -\frac1z + 2\beta(2z+1)\right) \mathrm{B}\left(z,\frac12\right) ,
\end{split}
\end{equation}
where $\psi(z)$ is the polygamma function and $\beta(z)=\frac12(\psi(\frac{z+1}{2})-\psi(\frac{z}{2}))$. For the second equality in \eqref{derry} we have used $\psi(z+1)=\psi(z)+\frac1z$. For $|z|<\frac12$ we have by \cite[8.373.1]{GR}
\begin{equation}\label{zero000}
-\frac1z + 2\beta(2z+1) = -\frac1z + 2 \log 2 +2 \sum_{k=1}^\infty (-1)^k (2^k-1) \zeta(k+1) z^k = \sum_{k=-1}^\infty c_k z^k ,
\end{equation}
where
\begin{equation}
c_{-1}=-1, \qquad c_0 = 2\log 2, \quad c_k =2(-1)^k (2^k-1) \zeta(k+1) \quad (k\geq 1).
\end{equation}
With this, \eqref{derry} gives 
\begin{equation}\label{bbcc}
(n+1) b_{n+1} =  \sum_{k=-1}^{n+1} c_{n-k} b_k \quad (n\geq -1) ,
\end{equation}
and we obtain the following recurrence relation for the Laurent coefficients of $\mathrm{B}\left(z,\frac12\right)$:
\begin{equation}\label{bbcc}
b_{n+1} = \frac{1}{n+2} \sum_{k=-1}^{n} c_{n-k} b_k \quad (n\geq 0) , \qquad \sum_{k=-1}^{0} c_{-k-1} b_k = 0 .
\end{equation}
The last relation yields $b_0=2\log 2$, and applying the recurrence we get
\begin{equation}
\begin{split}
b_1 & =2 (\log 2)^2-\frac{\pi ^2}{6} \\
b_2 & =2 \zeta (3)+ \frac43 (\log 2)^3-\frac{\pi ^2}{3} \log 2 \\
b_3 & = 4 \zeta (3) \log 2-\frac{\pi ^4}{40}+\frac{2}{3} (\log 2)^4-\frac{\pi ^2}{3} (\log 2)^2\\
b_4 & = 6 \zeta (5)+\left( 4 (\log 2)^2-\frac{\pi ^2}{3}\right) \zeta (3)+\frac{4 }{15} (\log 2)^5-\frac{2\pi ^2}{9}  (\log 2)^3-\frac{\pi ^4}{20}  \log 2,
\end{split}
\end{equation}
and so on.

As to the second term on the right-hand side of \eqref{MMalpha},
\begin{equation}
\frac{2^{-\alpha} - 1 + \alpha \log 2}{\alpha^2}  = \sum_{k=0}^\infty \frac{(-\log 2)^{k+2}\alpha^k}{(k+2)!}
\end{equation}
and so
\begin{equation}
\left. \frac{d^n}{d\alpha^n}  \frac{2^{-\alpha} - 1 + \alpha \log 2}{\alpha^2}  \right|_{\alpha=0}
= \frac{(-\log 2)^{n+2}}{(n+2)(n+1)} .
\end{equation}
We conclude
\begin{equation}\label{rhonny}
\rho_n = \frac{24}{\pi^2} \left( \frac{n!  (-1)^{n+1} b_{n+1}}{2^{n+1} } + \frac{2(-\log 2)^{n+2}}{(n+2)(n+1)}  \right).
\end{equation}

In view of \eqref{rhonny},
\begin{equation}
\begin{split}
\rho_1 & = \frac{12}{\pi^2} \zeta(3) -2\log 2 = 0.0752316\ldots \\
\rho_2 & = \frac{3\pi^2}{20} +2 (\log 2)^2 - \frac{24}{\pi^2} \zeta(3) \log2 =0.415242\ldots \\
\rho_3 & = -3 \zeta (3)+\frac{54}{\pi^2} \zeta (5)+\frac{36}{\pi^2} \zeta (3) (\log 2)^2-2 (\log 2)^3-\frac{9}{20} \pi ^2 \log 2 = 0.429262\ldots
.
\end{split}
\end{equation}

We now use the inversion formulas \eqref{inv0} and \eqref{inversion} (combined as in \eqref{inv0009}) to obtain
\begin{equation}
\mu_1 = \frac12 \rho_1 +\log 2 -\frac12 , \quad 
\mu_2 = \frac12 \rho_2 + \rho_1 (\log 2- \tfrac12) - \log 2 +(\log 2)^2 +\frac12 , \quad \text{etc.,}
\end{equation}
and the relations stated in \eqref{explicit}--\eqref{explicit3} follow.

\begin{figure}
\begin{center}
\includegraphics[width=0.435\textwidth]{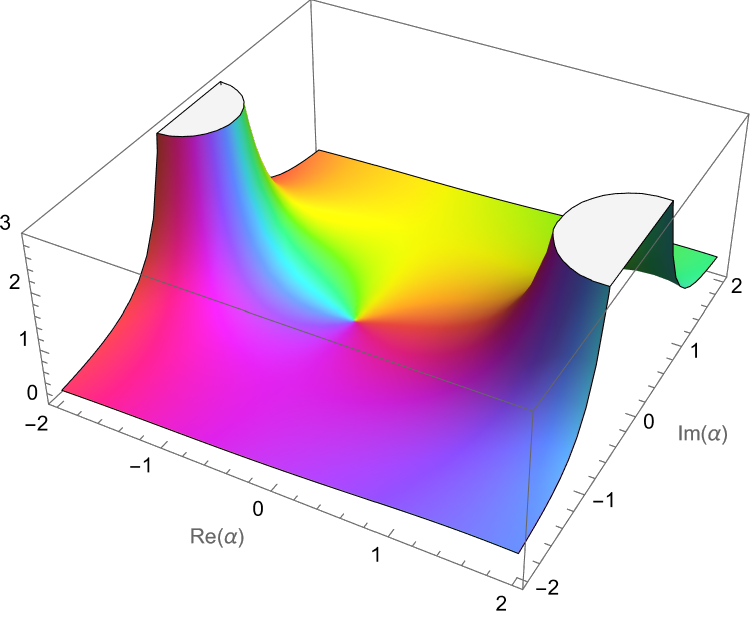}\quad\includegraphics[width=0.435\textwidth]{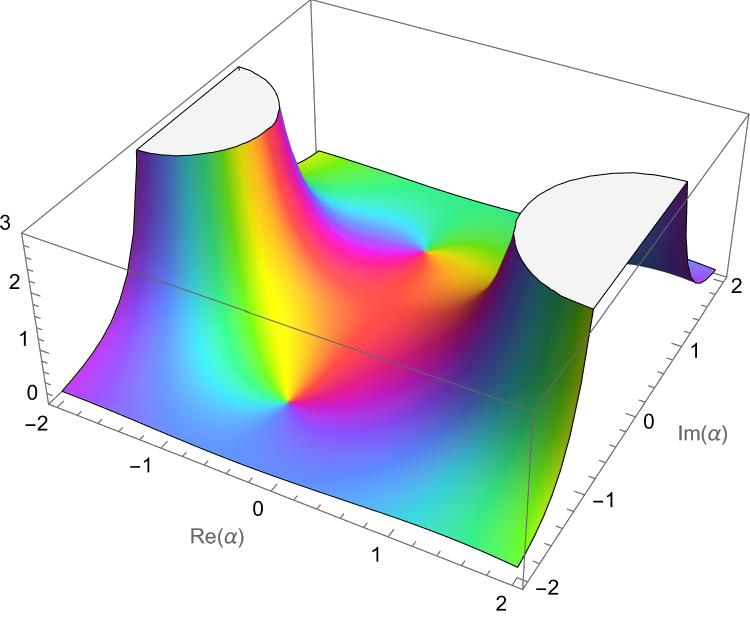}\quad\includegraphics[width=0.07\textwidth]{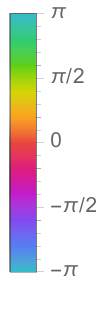}
\end{center}
\caption{The first and second derivatives of $M(\alpha)$, respectively, representing the generalized moments 
$\mu_{1,\alpha}$ and $\mu_{2,\alpha}$ in \eqref{genmu}. The height of the graph represents the function's absolute value and the colour its argument.} \label{fig2}
\end{figure}

\begin{sidewaystable}[!htbp]
\begin{center}
\vspace{400pt}
\begin{tabular}{|c|cccc|}
\hline
 $\alpha$ & $\mu_{0,\alpha}$ & $\mu_{1,\alpha}$ & $\mu_{2,\alpha}$ & $\mu_{3,\alpha}$   \\
\hline & & & &  \\[-10pt]
$-1$ & $\displaystyle \frac{12 (4-\pi)}{\pi ^2}$ & $\displaystyle \frac{12 (4-\pi  \log 4)}{\pi ^2}$ & $\displaystyle \frac{192-\pi  (\pi ^2+24+48 \log ^2 2)}{\pi ^2}$ & $\displaystyle \frac{6 (96-3 \pi  \zeta (3)-\pi  \log 2\, (\pi ^2+24+16 \log ^2 2))}{\pi ^2}$ \\[5pt]
$0$ & $1$ & $\displaystyle \frac{6 \zeta (3)}{\pi ^2}-\frac{1}{2}$ & $\displaystyle \frac{3 \pi ^2}{40}+\frac{1}{2}-\frac{6 \zeta (3)}{\pi ^2}$ & 
$\displaystyle \frac{9 (\zeta (3)+3 \zeta (5))}{\pi ^2}-\frac{9 \pi ^2}{80}-\frac{6 \zeta (3)+3}{4}$\\[5pt]
$1$ & $\displaystyle \frac{16}{\pi ^2}$ & $\displaystyle \frac{16 (3 \pi -7)}{3 \pi ^2}$ & $\displaystyle \frac{32 (34+3 \pi  (6\log 2-7))}{9 \pi ^2}$ & $\displaystyle \frac{4 (9 \pi ^3-1136+48 \pi  (17+9 \log ^2 2-21 \log 2))}{9 \pi ^2}$ \\[5pt]
\hline
\hline
 $\alpha$ & \multicolumn{2}{c}{$\mu_{0,\alpha}$} & \multicolumn{2}{c|}{$\mu_{2,\alpha}$}   \\
\hline &  &  & &  \\[-10pt]
$-\frac12$ & \multicolumn{2}{c}{$\frac{16}{\pi ^2} \left(8-\frac{\sqrt{2 \pi } \Gamma \left(\frac{1}{4}\right)}{\Gamma \left(\frac{3}{4}\right)}\right)$} 
& \multicolumn{2}{c|}{$\frac{16}{9 \pi ^2} \left(1408-\frac{\sqrt{2 \pi } \Gamma \left(\frac{5}{4}\right) \left(144 \mathbf{G} +9 \pi ^2+224+12 \pi  (4+3\log 2)+12 \log 2 (8+3\log 2)\right)}{\Gamma \left(\frac{3}{4}\right)}\right)$}  \\[10pt]
$\frac12$ & \multicolumn{2}{c}{$\frac{96}{5 \pi ^2} \left(4+\frac{\sqrt{2 \pi } \Gamma \left(-\frac{1}{4}\right)}{\Gamma \left(\frac{1}{4}\right)}\right)$} & \multicolumn{2}{c|}{$\frac{24}{125 \pi ^2} \left(11008+\frac{\sqrt{2 \pi } \Gamma \left(-\frac{1}{4}\right) \left(-400 \mathbf{G}+25 \pi ^2+440 \pi +2752+100 \log ^2 2-20 (44+5 \pi ) \log 2\right)}{\Gamma \left(\frac{1}{4}\right)}\right)$}  \\[10pt]
\hline
\end{tabular}
\end{center}
\caption{Examples of explicit values of the generalized moment $\mu_{n,\alpha}$ computed from \eqref{genmu} using Mathematica, where $\mathbf{G}$ denotes Catalan's constant. The standard moment corresponds to $n=0$ and the logarithmic moment to $\alpha=0$. \label{table1}}
\end{sidewaystable}

The same denominated convergence argument used for logarithmic moments also applies to more general test functions; again a direct corollary of \cite[Proposition 2]{Smallest}. This shows that, for any interval $\scrD\subset[0,1]$, $0\leq a<2$, $C>0$ and $F:\RR_{>0}\to\CC$ continuous with $|F(s)| \leq   C\max(s^a,s^{-a})$, we have
\begin{equation}\label{LD1MomMD}
\lim_{\delta\to 0}   \int_{\scrD}   F\left(\delta^{1/2} q_{\min}(x,\delta)\right) dx=
\vol\scrD \int_0^\infty F(s) \eta(s)\, ds.
\end{equation}

The same holds in the case of discrete sampling, and also in arbitrary dimension $d$ where $0\leq a<d+1$, as corollaries of the results in \cite[Sections 2 and 3]{Smallest}.

Admissible test functions include the generalized moments $F(s)=s^\alpha (\log s)^n $ with $n$ any non-negative integer and $|\Re\alpha|<2$. This generalized moment gives the standard moment for $n=0$, and the logarithmic moment for $\alpha=0$. We have for the limit
\begin{equation}\label{genmu}
\mu_{n,\alpha} =\int_0^\infty s^\alpha (\log s)^n \eta(s)\, ds = \frac{d^n}{d\alpha^n} \int_0^\infty s^\alpha \eta(s)\, ds =  \frac{d^n M(\alpha)}{d\alpha^n} .
\end{equation}
The graphs of these functions are displayed in Fig.~\ref{fig2} for $n=1$ and $2$; for $n=0$ see \cite[Fig.~2]{Smallest}. Table \ref{table1} comprises examples of explicit values of $\mu_{n,\alpha}$ computed from \eqref{genmu} using Mathematica. The case $\alpha=0$ corresponds to the logarithmic moments $\mu_n$ calculated ``by hand'' in the previous section.

\section{Minimal Resonance Orders in Higher Dimensions}\label{secMS}

Let us now turn to questions posed in the conclusions of \cite[Section 7]{Meiss23}. Following \cite{Meiss21,Meiss23} we define the minimal resonance order of a vector $\vecomega\in\RR^d$ by
\begin{equation}\label{MRO}
M(\vecomega,\delta) = \min_{\vecp\in\ZZ^d\setminus\{\vecnull\}} \left\{ \| \vecp \|_1 :  \min_{q\in\ZZ} \Delta_{\vecp,q}(\vecomega) \leq \delta \right\},
\end{equation}
where
\begin{equation}
\Delta_{\vecp,q}(\omega) =\frac{|\vecp\cdot\vecomega - q|}{\| \vecp\|_2} .
\end{equation}
The quantity $M(\vecomega,\delta)$ is a measure of how close $\vecomega$ is to a rational vector and arises naturally in the characterization of breakdowns of invariant tori in integrable systems under perturbation. In dimension one this reduces to $q_{\min}(\omega,2\delta)$.

For $L>0$, let 
\begin{equation}
R(L) = \mu\{ g\in\GamG : (\ZZ^{d+1}\setminus\{\vecnull\})\, g \cap \fB_L  = \emptyset\} ,
\end{equation}
where $G=\SL(d+1,\RR)$, $\Gamma=\SL(d+1,\ZZ)$, and
\begin{equation}
\fB_L=\{ (\vecx,y)\in \RR^{d+1} : \|\vecx\|_1\leq L, \; |y|\leq \|\vecx\|_2  \} .
\end{equation}
In \cite{Smallest} we defined $R(L)$ using primitive lattice points, but given the shape of $\fB_L$ both definitions are equivalent.

We will state the limit law for a general class of test functions $F$. We obtain the asymptotic value distribution as in \cite[Proposition 12]{Smallest} if we take $F$ to be the the indicator function of $\RR_{>L}$, and the generalized moments for $F(L)=L^\alpha (\log L)^n$ with $|\Re\alpha|<d+1$, $n=0,1,2,\ldots$.

\begin{proposition}
Let $\scrD\subset[0,1]^d$ with boundary of Lebesgue measure, $0\leq a <d+1$, $C>0$ and $F:\RR_{\geq 0}\to\CC$ continuous such that
\begin{equation}\label{unbounded007}
|F(L)| \leq C \max(L^{a},L^{-a}).
\end{equation}
Then
\begin{equation}\label{eqrmom}
\lim_{\delta\to 0}    \frac{1}{\vol\scrD} \int_\scrD  F\left( \delta^{1/(d+1)} M(\vecomega,\delta)\right) \,d\vecomega
= \int_0^\infty F(L) \;  dR(L) .
\end{equation}
\end{proposition}

\begin{proof} We first assume $F$ is the indicator function of an interval, and follow the same steps as in \cite[Proposition 12]{Smallest}. We note that
\begin{equation} \label{cal100011}
\begin{split}
&M(\vecomega,\delta) > L \delta^{-1/(d+1)}\\
&  \Leftrightarrow
 \left\{ (\vecp,q)\in\ZZ^{d+1} :  0<\| \vecp \|_1 \leq L \delta^{-1/(d+1)},\; \Delta_{\vecp,q}(\omega) \leq \delta \right\} =\emptyset \\
 &  \Leftrightarrow 
 \left\{ (\vecp,q)\in\ZZ^{d+1}\setminus\{ \vecnull \} :  \| \vecp \|_1 \leq L \delta^{-1/(d+1)},\; |\vecp\cdot\vecomega - q|\leq \delta \| \vecp\|_2  \right\} =\emptyset .
  \end{split}
\end{equation}
The last statement is in turn equivalent to
\begin{equation}\label{lala}
(\ZZ^{d+1}\setminus\{\vecnull\})  \cap \fB_L \begin{pmatrix} \delta^{-1/(d+1)} 1_d & \trans\vecnull \\ \vecnull & \delta^{d/(d+1)} \end{pmatrix} 
\begin{pmatrix} 1_d & \trans\vecomega \\ \vecnull & 1\end{pmatrix} = \emptyset.
\end{equation}
We can now directly apply \cite[Theorem 6.5, $\vecalf=\vecnull$]{partI} (cf.~also \cite{Artiles23}) to \eqref{lala}, which yields statement when $F$ is a indicator function of an interval. The extension to general bounded continuous functions follows by a standard approximation argument using finite linear combinations of indicator functions. 

Turning now to unbounded $F$, first of all note that $\fB_L=L\fB_1$. This means the lattice 
\begin{equation}\label{lalalalalala2}
\ZZ^{d+1} \begin{pmatrix} 1_d & -\trans\vecomega \\ \vecnull & 1\end{pmatrix} \begin{pmatrix} \delta^{1/(d+1)} 1_d & \trans\vecnull \\ \vecnull & \delta^{-d/(d+1)}
 \end{pmatrix} 
\end{equation}
needs to avoid the ball $L\scrB_0$, where $\scrB_0$ is any choice of open ball contained in $\fB_1$ and not containing the origin. By the argument in the proof of Proposition 4, after (2.22) in \cite{Smallest} (replacing all matrices $g$ by $\trans g^{-1}$), we obtain the estimate, valid for all $L\geq 1$, $0<\delta\leq 1$ and some sufficiently large constant $C$,
 \begin{equation}\label{esc1}
\vol\left\{ \vecomega\in\scrD : \delta^{1/(d+1)}  M(\vecomega,\delta) > L \right\} \leq C L^{-(d+1)} .
\end{equation}
A key input here, as explained in \cite{Smallest}, is the escape-of-mass estimate provided in \cite{Kim23}.
The next step is to control small values of $L$, and thus the measure of $\vecomega$ for which $L\fB_1$ contains a non-zero element of the lattice \eqref{lalalalalala2}. Hence the shortest non-zero vector of \eqref{lalalalalala2} has to have length $\ll L$, and following the relevant steps of \cite{Smallest} in the proof of its Proposition 4, we obtain 
\begin{equation}\label{esc2}
\vol\left\{ \vecomega\in\scrD : \delta^{1/(d+1)}  M(\vecomega,\delta) \leq L \right\} \leq C L^{d+1} ,
\end{equation}
valid for all $0<L\leq 1$, $0<\delta\leq 1$ and a constant $C$; cf.~\cite[(2.31)]{Smallest}. The estimates \eqref{esc1} and \eqref{esc2} now permit the extension of \eqref{eqrmom} to test functions $F$ subject to \eqref{unbounded007}.
\end{proof}

The density of $R(L)$ corresponds to the histogram in Fig.~8(a) of \cite{Meiss21} for $d=2$. The exponent of $\delta$  in Fig. 8(b) is approximately $0.334$, which is consistent with the theoretically predicted scaling with exponent $1/(d+1)$. Note furthermore that the choice $F(L)=\log L$ gives the small $\delta$ asymptotics of the expectation
\begin{equation}\label{eqrmom22}
\frac{1}{\vol\scrD} \int_\scrD \log M(\vecomega,\delta)  \,d\vecomega
= -\frac{\log\delta}{d+1} + \int_0^\infty \log L \; dR(L) +o(1) ,
\end{equation}
which is compatible with the numerics of \cite{Meiss23} in dimension $d=2$.

\vskip20pt\noindent {\bf Acknowledgements.}
I thank the referee for helpful comments and suggestions. This research was supported by EPSRC grant EP/S024948/1. Data supporting this study are included within the article.

\end{document}